\documentclass[oneside,english]{amsart}
\usepackage[T1]{fontenc}
\usepackage[latin9]{inputenc}
\usepackage{amstext}
\usepackage{amsthm}
\usepackage{amssymb}

\makeatletter

\providecommand{\tabularnewline}{\\}

\numberwithin{equation}{section}
\numberwithin{figure}{section}
\theoremstyle{plain}
\newtheorem{thm}{\protect\theoremname}

\makeatother

\usepackage{babel}
\providecommand{\theoremname}{Theorem}

\begin{document}
\title{EULER'S AND THE TAXI CAB RELATIONS AND OTHER NUMBERS THAT CAN BE WRITTEN
TWICE AS SUMS OF TWO CUBED INTEGERS}
\author{Vladimir PLETSER}
\address{European Space Agency (ret.)}
\email{Pletservladimir@gmail.com}
\begin{abstract}
We show that Euler's relation and the Taxi-Cab relation are both solutions
of the same equation. General solutions of sums of two consecutive
cubes equaling the sum of two other cubes are calculated. There is
an infinite number of relations to be found among the sums of two
consecutive cubes and the sum of two other cubes, in the form of two
families. Their recursive and parametric equations are calculated.
\end{abstract}


\maketitle

\section{Introduction\label{sec:Introduction}}

The remarkable relation 
\begin{equation}
3^{3}+4^{3}+5^{3}=6^{3}\label{eq:1}
\end{equation}
among the cubes of four successive integers is often attributed to
Euler, while in fact it was already known to P. Bungus in the XVIth
century \cite{Bungus,Dickson}. No other similar relation can be found
between cubes of four successive integers. 

Another well-known relation involving two different sums of two cubes
is
\begin{equation}
1729=9^{3}+10^{3}=12^{3}+1^{3}\label{eq:2}
\end{equation}
often call the taxi-cab number or taxi-cab relation and attributed
to Indian mathematician Ramanujan after he mentioned in 1919 to fellow
British mathematician Hardy that this number is remarkable in the
fact that it is the smallest integer that can be written as the sum
of two positive cubes in two ways (see historical account in e.g.,
\cite{Grinstein}). However, this relation was already mentioned by
French mathematician Frenicle in the XVIIth century \cite{Dickson,Frenicle}.
Nevertheless, we will refer in this paper to (\ref{eq:1}) and (\ref{eq:2})
as Euler's relation and Ramanujan's taxi-cab relation.

In fact, both relations can be deduced from a same equation, as we
show in this paper. It is simple to see that one can find other taxi-cab
numbers smaller than Ramanujan's by transferring one term from left
to right of (\ref{eq:1}), introducing negative integers and yielding
successively
\begin{align}
91 & =3^{3}+4^{3}=6^{3}+\left(-5\right)^{3}\label{eq:3-1}\\
152 & =3^{3}+5^{3}=6^{3}+\left(-4\right)^{3}\label{eq:4-1}\\
189 & =4^{3}+5^{3}=6^{3}+\left(-3\right)^{3}\label{eq:5-1}
\end{align}
and so on. Other taxi-cab numbers can be found by multiplying each
relations (\ref{eq:3-1}) to (\ref{eq:5-1}) by $k^{3}$, i.e., the
cube of any integer $k$. Sequences A001235 and A051347 in the OEIS
\cite{OEIS} list all taxi-cab numbers for respectively only positive
integers and for positive and negative integers. Numerous mathematicians
and authors have worked on sums of cubes and equal sums of cubes.
Excellent summaries and numerous results can be found e.g., in \cite{Dickson,Piezas}.

In this paper, our interest is in numbers that can be written as sums
of two cubes in at least two ways, one of them involving two consecutive
cubes. In Section, 2 we show first that Euler's and the Taxi-Cab relations
are solutions of the same equation. We calculate then the general
case of the sum of two consecutive cubes equal to the sum of two other
cubes. In Section 3, we characterize two infinite families of solutions
of the sum of two consecutive cubes equaling the sum of two other
cubes.

\section{General equation\label{sec:General-equation}}

We show first that (\ref{eq:1}) and (\ref{eq:2}) are both solutions
of a same equation. If one observes that the first term on the right
hand side in (\ref{eq:2}) and (\ref{eq:3-1}) is three units larger
than the first term on the left hand side, we can write
\begin{equation}
N=n^{3}+\left(n+1\right)^{3}=\left(n+3\right)^{3}+\left(n+\alpha\right)^{3}\label{eq:10}
\end{equation}
with $\alpha$ integer and $N$ the positive integer that can be represented
in (at least) two ways by sums of two consecutive cubes and of two
other cubes, one of which possibly negative. Equation (\ref{eq:10})
yields the two general solutions 
\begin{equation}
n=\frac{-3\left(\alpha^{2}+8\right)\pm\sqrt{-3\left(\alpha^{4}+8\left(\alpha^{3}-6\alpha^{2}+13\alpha+2\right)\right)}}{6\left(\alpha+2\right)}\label{eq:11}
\end{equation}
which produces integer solutions for $\alpha=-8$, giving $n_{+}=3$
and $n_{-}=9$. Equation (\ref{eq:10}) yields then respectively (\ref{eq:3-1})
and (\ref{eq:2}), showing that (\ref{eq:10}) yields both Euler's
relation and Ramanujan's taxi-cab number relation.

Let us consider now the general equation
\begin{equation}
N=n^{3}+\left(n+1\right)^{3}=\left(n+a\right)^{3}+\left(n+b\right)^{3}\label{eq:3}
\end{equation}
with $a$ and $b$ integers, $a>0$ and $b<0$. Solving for $n$,
the third degree equation defined by the second equality in (\ref{eq:3})
reduces to a second degree equation 
\begin{equation}
3n^{2}\left(a+b-1\right)+3n\left(a^{2}+b^{2}-1\right)+\left(a^{3}+b^{3}-1\right)=0\label{eq:4}
\end{equation}
whose discriminant reads
\begin{equation}
D=3\left(\left(a-b\right)^{4}-\left(a^{4}+b^{4}+\left(a-1\right)^{4}+\left(b-1\right)^{4}\right)+1\right)\label{eq:5}
\end{equation}
providing two real solutions for $D>0$, namely

\begin{equation}
n=\frac{-3\left(a^{2}+b^{2}-1\right)\pm\sqrt{3\left(\left(a-b\right)^{4}-\left(a^{4}+b^{4}+\left(a-1\right)^{4}+\left(b-1\right)^{4}\right)+1\right)}}{6\left(a+b-1\right)}\label{eq:6}
\end{equation}
As $N$ must be positive, we limit our search to $b<0<a<\left|b\right|$
and discrete solutions of (\ref{eq:3}) or (\ref{eq:4}) are found,
as shown in Table \ref{tab:Table 1} for $n<1000$, and arranged in
increasing order of $N$. 

\begin{table}
\caption{\label{tab:Table 1}Values of $a$, $b$, $n_{+}$, $n_{-}$, solutions
of (\ref{eq:3}) for $n<1000$ and $b<0<a<\left|b\right|$}

\centering{}%
\begin{tabular}{|c|c|c|c|l|}
\hline 
$a$ & $b$ & $n_{+}$ & $n_{-}$ & $N=n^{3}+\left(n+1\right)^{3}=\left(n+a\right)^{3}+\left(n+b\right)^{3}$\tabularnewline
\hline 
\hline 
3 & -8 & 3 & -- & $91=3^{3}+4^{3}=6^{3}+(-5)^{3}$\tabularnewline
\hline 
2 & -7 & 4 & -- & $189=4^{3}+5^{3}=6^{3}+(-3)^{3}$\tabularnewline
\hline 
3 & -8 & -- & 9 & $1729=9^{3}+10^{3}=12^{3}+1^{3}$\tabularnewline
\hline 
10 & -39 & 18 & -- & $12691=18^{3}+19^{3}=28^{3}+(-21)^{3}$\tabularnewline
\hline 
9 & -38 & -- & 32 & $68705=32^{3}+33^{3}=41^{3}+(-6)^{3}$\tabularnewline
\hline 
10 & -39 & -- & 36 & $97309=36^{3}+37^{3}=46^{3}+(-3)^{3}$\tabularnewline
\hline 
105 & -194 & 46 & -- & $201159=46^{3}+47^{3}=151^{3}+(-148)^{3}$\tabularnewline
\hline 
32 & -127 & 58 & -- & $400491=58^{3}+59^{3}=90^{3}+(-69)^{3}$\tabularnewline
\hline 
64 & -243 & 107 & -- & $2484755=107^{3}+108^{3}=171^{3}+(-136)^{3}$ \tabularnewline
\hline 
73 & -258 & 108 & -- & $2554741=108^{3}+109^{3}=181^{3}+(-150)^{3}$ \tabularnewline
\hline 
32 & -103 & -- & 121 & $3587409=121^{3}+122^{3}=153^{3}+18^{3}$\tabularnewline
\hline 
248 & -481 & 121 & -- & $3587409=121^{3}+122^{3}=369^{3}+(-360)^{3}$ \tabularnewline
\hline 
37 & -192 & -- & 123 & $3767491=123^{3}+124^{3}=160^{3}+(-69)^{3}$\tabularnewline
\hline 
43 & -168 & -- & 163 & $8741691=163^{3}+164^{3}=206^{3}+(-5)^{3}$\tabularnewline
\hline 
91 & -360 & 163 & -- & $8741691=163^{3}+164^{3}=254^{3}+(-197)^{3}$ \tabularnewline
\hline 
819 & -1208 & 197 & -- & $15407765=197^{3}+198^{3}=1016^{3}+(-1011)^{3}$ \tabularnewline
\hline 
57 & -128 & -- & 235 & $26122131=235^{3}+236^{3}=292^{3}+107^{3}$\tabularnewline
\hline 
184 & -597 & 235 & -- & $26122131=235^{3}+236^{3}=419^{3}+(-362)^{3}$\tabularnewline
\hline 
77 & -208 & -- & 301 & $54814509=301^{3}+302^{3}=378^{3}+93^{3}$\tabularnewline
\hline 
120 & -629 & 393 & -- & $121861441=393^{3}+394^{3}=513^{3}+(-236)^{3}$ \tabularnewline
\hline 
120 & -629 & -- & 411 & $139361059=411^{3}+412^{3}=531^{3}+(-218)^{3}$\tabularnewline
\hline 
393 & -1178  & 438 & -- & $168632191=438^{3}+439^{3}=831^{3}+(-740)^{3}$\tabularnewline
\hline 
152 & -793 & 481 & -- & $223264809=481^{3}+482^{3}=633^{3}+(-312)^{3}$ \tabularnewline
\hline 
128 & -511 & -- & 490 & $236019771=490^{3}+491^{3}=618^{3}+(-21)^{3}$\tabularnewline
\hline 
3225 & -4274 & 528 & -- & $295233841=528^{3}+529^{3}=3753^{3}+(-3746)^{3}$\tabularnewline
\hline 
148 & -687 & -- & 562 & $355957875=562^{3}+563^{3}=710^{3}+(-125)^{3}$ \tabularnewline
\hline 
2258  & -3367 & 562 & -- & $355957875=562^{3}+563^{3}=2820^{3}+(-2805)^{3}$\tabularnewline
\hline 
512  & -1591 & 607 & -- & $448404255=607^{3}+608^{3}=1119^{3}+(-984)^{3}$\tabularnewline
\hline 
777  & -1952 & 633 & -- & $508476241=633^{3}+634^{3}=1410^{3}+(-1319)^{3}$\tabularnewline
\hline 
190 & -999 & -- & 640 & $525518721=640^{3}+641^{3}=830^{3}+(-359)^{3}$ \tabularnewline
\hline 
442  & -1767 & 804 & -- & $1041378589=804^{3}+805^{3}=1246^{3}+(-963)^{3}$\tabularnewline
\hline 
\end{tabular}
\end{table}

Note also that similar relations but with coefficients having opposite
signs are obtained for negative values of $a$ and for $a^{\prime}=-b+1$,
$b^{\prime}=-a+1$, and $n^{\prime}=-n-1$.

It is seen also that three sums of two cubes are found for $n=121,163,235,562.$
Other relations are given in OEIS \cite{OEIS} Sequences A352133 to
A352136 and cases with three sums of two cubes are given in Sequences
A352220 to A352225.

\section{Two Infinite Families\label{sec:Two-Infinite-Families}}

Figure 3.1 of \cite{Pletser} shows a plot of the couples $\left(n,n+a\right)$
for $0<n\leq275$ (data are from OEIS \cite{OEIS} Sequences A352135,
A352136, A352222, A352223, A352224, A352225). Two families are clearly
visible along two curves.

The first top curve (or first family) includes all couples $\left(n,n+a\right)$
such that $\eta=\left(n+a\right)+\left(n+b\right)=2n+a+b$ are regularly
increasing odd integers as shown in Table \ref{tab:Table 2} for the
first twenty cases, while for the second below curve (or second family),
$\eta=2n+a+b$ are regularly increasing odd multiples of $3$.

\begin{table}
\caption{\label{tab:Table 2}Values of $n$, $n+a$, $n+b$, $\eta=2n+a+b$
for first and second families}

\centering{}%
\begin{tabular}{|c||c|c|c|c||c|c|c|c|}
\hline 
 & \multicolumn{4}{c||}{First family} & \multicolumn{4}{c|}{Second family}\tabularnewline
\cline{2-9} \cline{3-9} \cline{4-9} \cline{5-9} \cline{6-9} \cline{7-9} \cline{8-9} \cline{9-9} 
$i$ & $n$ & $n+a$ & $n+b$ & $\eta$ & $n$ & $n+a$ & $n+b$ & $\eta$\tabularnewline
\hline 
\hline 
1 & 3 & 6 & -5 & 1 & 4 & 6 & -3 & 3\tabularnewline
\hline 
2 & 46 & 151 & -148 & 3 & 121 & 369 & -360 & 9\tabularnewline
\hline 
3 & 197 & 1016 & -1011 & 5 & 562 & 2820 & -2805 & 15\tabularnewline
\hline 
4 & 528 & 3753 & -3746 & 7 & 1543 & 10815 & -10794 & 21\tabularnewline
\hline 
5 & 1111 & 10090 & -10081 & 9 & 3280 & 29538 & -29511 & 27\tabularnewline
\hline 
6 & 2018 & 22331 & -22320 & 11 & 5989 & 65901 & -65868 & 33\tabularnewline
\hline 
7 & 3321 & 43356 & -43343 & 13 & 9886 & 128544 & -128505 & 39\tabularnewline
\hline 
8 & 5092 & 76621 & -76606 & 15 & 15187 & 227835 & -227790 & 45\tabularnewline
\hline 
9 & 7403 & 126158 & -126141 & 17 & 22108 & 375870 & -375819 & 51\tabularnewline
\hline 
10 & 10326 & 196575 & -196556 & 19 & 30865 & 586473 & -586416 & 57\tabularnewline
\hline 
11 & 13933 & 293056 & -293035 & 21 & 41674 & 875196 & -875133 & 63\tabularnewline
\hline 
12 & 18296 & 421361 & -421338 & 23 & 54751 & 1259319 & -1259250 & 69\tabularnewline
\hline 
13 & 23487 & 587826 & -587801 & 25 & 70312 & 1757850 & -1757775 & 75\tabularnewline
\hline 
14 & 29578 & 799363 & -799336 & 27 & 88573 & 2391525 & -2391444 & 81\tabularnewline
\hline 
15 & 36641 & 1063460 & -1063431 & 29 & 109750 & 3182808 & -3182721 & 87\tabularnewline
\hline 
16 & 44748 & 1388181 & -1388150 & 31 & 134059 & 4155891 & -4155798 & 93\tabularnewline
\hline 
17 & 53971 & 1782166 & -1782133 & 33 & 161716 & 5336694 & -5336595 & 99\tabularnewline
\hline 
18 & 64382 & 2254631 & -2254596 & 35 & 192937 & 6752865 & -6752760 & 105\tabularnewline
\hline 
19 & 76053 & 2815368 & -2815331 & 37 & 227938 & 8433780 & -8433669 & 111\tabularnewline
\hline 
20 & 89056 & 3474745 & -3474706 & 39 & 266935 & 10410543 & -10410426 & 117\tabularnewline
\hline 
\end{tabular}
\end{table}

\subsection{Recursive relations\label{subsec:Recursive-relations}}

The values of $n$, $n+a$ and $n+b$ of both the first and second
families can be found by the recurrence relations
\begin{align}
n_{i} & =3n_{i-1}-3n_{i-2}+n_{i-3}+\kappa\label{eq:12}\\
\left(n+a\right)_{i} & =3\left(n+a\right)_{i-1}-3\left(n+a\right)_{i-2}+\left(n+a\right)_{i-3}+\lambda\label{eq:12-1}\\
\left(n+b\right)_{i} & =3\left(n+b\right)_{i-1}-3\left(n+b\right)_{i-2}+\left(n+b\right)_{i-3}-\lambda\label{eq:12-2}
\end{align}
with $\kappa=72$ and $216$ and $\lambda=576\left(i-2\right)$ and
$1728\left(i-2\right)$ for respectively the first and second families,
and the first three values of $n_{i}$, $\left(n+a\right)_{i}$ and
$\left(n+b\right)_{i}$ from Table 2.

\subsection{Parametric relations\label{subsec:Parametric-relations}}

We see from Table 2 that the fourth term $n+b$ of (\ref{eq:3}) is
negative and is decreasing regularly with increasing $n$. So, let
us pose $n+b=-\left(n+a\right)+\beta$, yielding from (\ref{eq:3})
\begin{equation}
N=n^{3}+\left(n+1\right)^{3}=\left(n+a\right)^{3}-\left(n+a-\beta\right)^{3}\label{eq:13}
\end{equation}
For specific relations between $a$ and $n$, one obtains two infinite
families of solutions as shown in the following two theorems, giving
parametric solutions for $n$, $N$, $n+a$ and $n+b$.
\begin{thm}
For $\forall i\in\mathbb{Z}_{0}^{+}$, $\exists n,a,\beta\in\mathbb{Z}_{0}^{+}$,
such that
\begin{equation}
a=\left(\beta-1\right)n+\beta^{2}+\beta+1\label{eq:15-1}
\end{equation}
and an infinite family of solutions of (\ref{eq:13}) exists for $\beta$
odd, 
\begin{equation}
\beta=2i-1\label{eq:16-1}
\end{equation}
yielding
\begin{align}
n & =\frac{\left(2i-1\right)\left(3\left(2i-1\right)^{2}+4\right)-1}{2}\label{eq:17-1}\\
N & =\frac{\left(2i-1\right)\left(3\left(2i-1\right)^{2}+4\right)\left(\left(2i-1\right)^{2}\left(3\left(2i-1\right)^{2}+4\right)^{2}+3\right)}{4}\label{eq:18-1}\\
n+a & =\frac{3\left(2i-1\right)^{2}\left(\left(2i-1\right)^{2}+2\right)+2i+1}{2}\label{eq:18-2}\\
n+b & =-\frac{3\left(2i-1\right)^{2}\left(\left(2i-1\right)^{2}+2\right)-2i+3}{2}\label{eq:18-3}
\end{align}
\end{thm}

\begin{proof}
Let $n,a,\beta,i\in\mathbb{Z}_{0}^{+}$, and let $a$, $n$ and $\beta$
satisfy (\ref{eq:15-1}). Relation (\ref{eq:13}) yields then the
third degree equation
\begin{equation}
n^{3}+\left(n+1\right)^{3}-\left(\beta n+\beta^{2}+\beta+1\right)^{3}+\left(\beta n+\beta^{2}+1\right)^{3}=0
\end{equation}
that simplifies immediately in the product of a linear and a quadratic
relations
\begin{equation}
\left(2n-\beta\left(3\beta^{2}+4\right)+1\right)\left(n^{2}+\left(2\beta+1\right)n+\beta^{2}+\beta+1\right)=0
\end{equation}
As the discriminant of the right quadratic polynomial is always negative,
the quadratic equation yields two complex solutions of no interest
here. The right linear equation yield the only real solution
\begin{equation}
n=\frac{\beta\left(3\beta^{2}+4\right)+1}{2}
\end{equation}
As $n$ must be integer, $\beta$ cannot be even and must be odd,
$\beta=2i-1$, yielding (\ref{eq:17-1}) to (\ref{eq:18-3}).
\end{proof}
\begin{thm}
For $\forall i\in\mathbb{Z}_{0}^{+}$, $\exists n,a,\beta\in\mathbb{Z}_{0}^{+}$,
such that
\begin{equation}
a=\frac{\left(\beta-3\right)n+2\beta}{3}\label{eq:15-1-1}
\end{equation}
and an infinite family of solutions of (\ref{eq:13}) exists for $\beta\equiv3\text{mod\ensuremath{6}}$,
\begin{equation}
\beta=3\left(2i-1\right)\label{eq:16-1-1}
\end{equation}
yielding
\begin{align}
n & =\frac{9\left(2i-1\right)^{3}-1}{2}\label{eq:17-1-2}\\
N & =\frac{27\left(2i-1\right)^{3}\left(27\left(2i-1\right)^{6}+1\right)}{4}\label{eq:17-2}\\
n+a & =\frac{3\left(2i-1\right)\left(3\left(2i-1\right)^{3}+1\right)}{2}\label{eq:17-3}\\
n+b & =-\frac{3\left(2i-1\right)\left(3\left(2i-1\right)^{3}-1\right)}{2}\label{eq:17-4}
\end{align}
\end{thm}

\begin{proof}
Let $n,a,\beta,i\in\mathbb{Z}_{0}^{+}$, and let $a$, $n$ and $\beta$
satisfy (\ref{eq:15-1-1}). Relation (\ref{eq:13}) yields then the
third degree equation
\begin{equation}
n^{3}+\left(n+1\right)^{3}-\left(\frac{\beta}{3}\left(n+2\right)\right)^{3}+\left(\frac{\beta}{3}\left(n-1\right)\right)^{3}=0\label{eq:18}
\end{equation}
that simplifies immediately in the product of a linear and a quadratic
relations
\begin{equation}
\frac{\left(6n-\beta^{3}+3\right)\left(n^{2}+n+1\right)}{3}=0\label{eq:19}
\end{equation}
The right quadratic equation yields two complex solutions of no interest
here. The right linear equation yield the only real solution
\begin{equation}
n=\frac{\beta^{3}-3}{6}\label{eq:20}
\end{equation}
As $n$ must be integer, $\beta$ must $3\text{mod\ensuremath{6}}$,
$\beta=3\left(2i-1\right)$, yielding (\ref{eq:17-1-2}) to (\ref{eq:17-4}).
\end{proof}

\section{Conclusion\label{sec:Conclusion}}

We have shown first that Euler's relation and the Taxi-Cab relation
are both solutions of the same equation. We have then calculated general
solutions of sums of two consecutive cubes equaling the sum of two
other cubes. We have finally shown that there is an infinite number
of relations that can be found among the sums of two consecutive cubes
and the sum of two other cubes, in the form of two families and we
have given their recursive and parametric equations.

\section*{Acknowledgment\label{sec:Acknowledgment}}

The Author wishes to acknowledge the help of an OEIS Associate Editor
and Editor-in-Chief, for additional computing in OEIS \cite{OEIS}
Sequence A352135.


\begin{thebibliography}{1}
\bibitem{Bungus}Bungus P. (1591). \emph{Numerorum Mysteria, 1618},
463; Pars Altera, 65.

\bibitem{Dickson}Dickson L.E. (2005). \emph{History of the Theory
of Numbers, Vol. II: Diophantine Analysis}, Dover Publications, New
York, 550-562.

\bibitem{Grinstein}Grinstein A. ( 2022). \emph{Ramanujan and 1729},
University of Melbourne Dept. of Math and Statistics Newsletter: Issue
3, 1998. available at https://web.archive.org/web/20040320144821/http://zadok.org/mattandloraine/1729.html,
Last accessed 17 March 2022.

\bibitem{Frenicle}Frenicle de Bessy B.(1657). \emph{Commercium Epistolicum
de Wallis, letter X}, Brouncker to Wallis, Oct. 13, 1657.

\bibitem{OEIS}Sloane N.J.A., ed. (2022). \emph{The On-Line Encyclopedia
of Integer Sequences}, published electronically at https://oeis.org.

\bibitem{Piezas}Piezas III T. (2010). \emph{A Collection of Algebraic
Identities, Chap 6: Third Powers}, available at https://sites.google.com/site/tpiezas/Home,
Last accessed 2 April 2022.

\bibitem{Pletser}Pletser V. (2022). \emph{Euler's and the Taxi-Cab
relations and other numbers that can be written twice as sums of two
cubed integers}, submitted. Preprint available at https://www.researchgate.net/publication/359706361\_EULER'S\_AND\_THE\_TAXI\_CAB\_RELATIONS\_AND\_OTHER\_NUMBERS\_THAT\_CAN\_BE\_WRITTEN\_TWICE\_AS\_SUMS\_OF\_TWO\_CUBED\_INTEGERS
\end{thebibliography}
\end{document}